\documentclass[a4paper]{article}
\usepackage{amsmath,amssymb}
\usepackage[latin1]{inputenc}
\usepackage[dvips]{epsfig}

\newcommand{\kk}{\mathbf{k}}

\newcommand{\Id}{\operatorname{Id}}
\newcommand{\zer}{\widehat{0}}
\newcommand{\one}{\widehat{1}}
\newcommand{\gch}{\epsfig{file=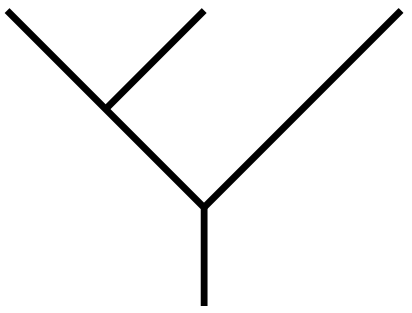,height=2mm}}
\newcommand{\drt}{\epsfig{file=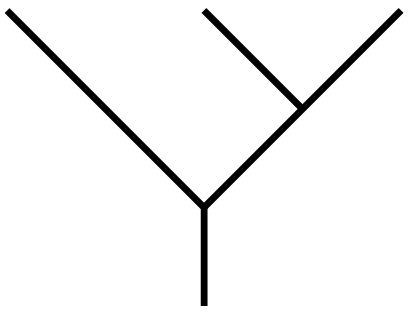,height=2mm}}

\newcommand{\dend}{\mathcal{Y}}
\newcommand{\YY}{\mathbb{Y}}
\newcommand{\TT}{\mathbb{T}}
\newcommand{\PP}{\mathcal{P}}
\newcommand{\Poset}{\mathbb{P}}
\renewcommand{\mod}{\operatorname{mod}}
\newcommand{\Dmod}{\mathcal{D}\mod}

\newtheorem{theorem}{Theorem}[section] 
\newtheorem{proposition}[theorem]{Proposition} 
 
\newtheorem{corollary}[theorem]{Corollary} 
\newtheorem{lemma}[theorem]{Lemma}

\newenvironment{proof}{\begin{trivlist}\item{\bf{Proof.}}}
  {\hfill\rule{2mm}{2mm}\end{trivlist}}

\title{On the Coxeter transformations \\ for Tamari posets }
\author{Frédéric Chapoton}
\date{\today}

\begin{document}

\maketitle

\begin{abstract}
  A relation between the anticyclic structure of the dendriform operad
  and the Coxeter transformations in the Grothendieck groups of the
  derived categories of modules over the Tamari posets is obtained.
\end{abstract}

\setcounter{section}{-1}

\section{Introduction}

There are now several algebraic structures on planar binary
trees. First, there is an operad, called the dendriform operad, whose
structure can be described by insertion of planar binary trees. Then
the free dendriform algebra on one generator is also an associative
algebra and in fact a Hopf algebra, called the Hopf algebra of planar
binary trees. Both the dendriform operad and the Hopf algebra of
planar binary trees have been shown to be related to a family of
posets on planar binary trees, called the Tamari lattices. 

Until recently, it was not realized that the dendriform operad is an
anticyclic operad. This fact implies the existence of a linear map of
order $n+1$ on the vector space spanned by planar binary trees with
$n+1$ leaves. The matrix of this endomorphism seemed similar to a
matrix appearing in the study of the Hopf algebra of planar binary
trees made in \cite{hivnov}. This was the starting point for this
article. 

The main result shows that the linear maps obtained from the
anticyclic structure of the dendriform operad can alternatively be
described using only the Tamari posets. More precisely, recall that,
for a quiver, the Coxeter transformation is the action induced on the
Grothendieck group by a canonical self-equivalence, called the
Auslander-Reiten translation, of the derived category of modules on
the quiver. Considering Tamari posets as quivers with relations gives
a family of Coxeter transformations on vector spaces spanned by planar
binary trees. Our result show that, up to sign, iterating twice the
Coxeter transformations recovers the anticyclic structure maps. All
this should hint at a deeper relationship between the dendriform
operad and derived categories for Tamari posets.

Also, this implies that the Coxeter transformation for Tamari posets
is periodic. It is expected that something similar should happen for
any Cambrian lattice associated to a finite Coxeter group
\cite{reading,thomas}.  More precisely, the Coxeter transformation in
the Grothendieck group of the derived category of modules on a
Cambrian lattice should have order dividing $2h+2$ where $h$ is the
Coxeter number of the Coxeter group. 

Let us also note that a similar, but much simpler and less interesting
theory can be done relating the diassociative anticyclic operad on one
hand and the family of total orders or chains on the other hand.

The article starts with many recollections on trees, posets, algebras,
operads and quivers. The main theorem and its proof are to be found in
section \ref{principale}. 

\textbf{Acknowledgements:} I would like to thank F. Hivert and J.-C. 
Novelli for stimulating discussions on the Hopf algebra of planar
binary trees. 

\section{Planar binary trees}

Let $n$ be a nonnegative integer. A planar binary tree of degree $n$
is a graph embedded in the plane which is a tree, has $n$ trivalent
vertices, $n+2$ univalent vertices and a distinguished univalent
vertex called the root. The other univalent vertices are called the
leaves. From now on, trivalent vertices and vertices will mean the
same thing. Planar binary trees are pictured with their root at the
bottom and leaves at the top, see Figure \ref{tam2}. 

Let $\YY(n)$ be the set of planar binary trees of degree $n$. It is a
classical combinatorial fact that the cardinality of $\YY(n)$ is the
Catalan number $c_n=\frac{1}{n+1}\binom{2n}{n}$. 

Let $\YY$ be the set of all planar binary trees and $\YY^+$ the set of
all planar binary trees but the tree $|$ with no vertex. For $S$ in
$\YY$, let $|S|$ be the degree of $S$, \textit{i.e.} its number of
vertices.  Let $Y$ be the unique tree with one vertex. 

Let us define some combinatorial operations on $\YY$. Let $S$ and $T$
be in $\YY$. Then let $S_{Y} T$ be the planar binary tree obtained by
grafting simultaneously $S$ to the left leaf of $Y$ and $T$ to the
right leaf of $Y$. This tree has degree $|S|+|T|+1$. 

Let $S \slash T$ be the tree obtained by grafting the root of $S$ to the
leftmost leaf of $T$. It has degree $|S|+|T|$. Similarly let $S
\backslash T$ be the tree obtained by grafting the root of $T$ to the
rightmost leaf of $S$. It also has degree $|S|+|T|$. 

Remark that one can also define $S_{Y} T$ as $(S/Y) \backslash T$ or
$S/ (Y \backslash T)$. The tree $|$ is a two-sided unit for both
$\backslash$ and $/$. 

There is an obvious involution on planar binary trees, given by the
left-right reversal of the plane. 
\section{Tamari posets}

There is a natural order relation $\leq$ on the set $\YY(n)$, which
was introduced and studied by Tamari in \cite{tamari}. 

The order relation $\leq$ is defined as the transitive closure of some
covering relations. A tree $S$ is covered by a tree $T$ if they differ
only in some neighborhood of an edge by the replacement of the
configuration $\gch$ in $S$ by the configuration $\drt$ in $T$.

\begin{figure}
  \begin{center}
    \epsfig{file=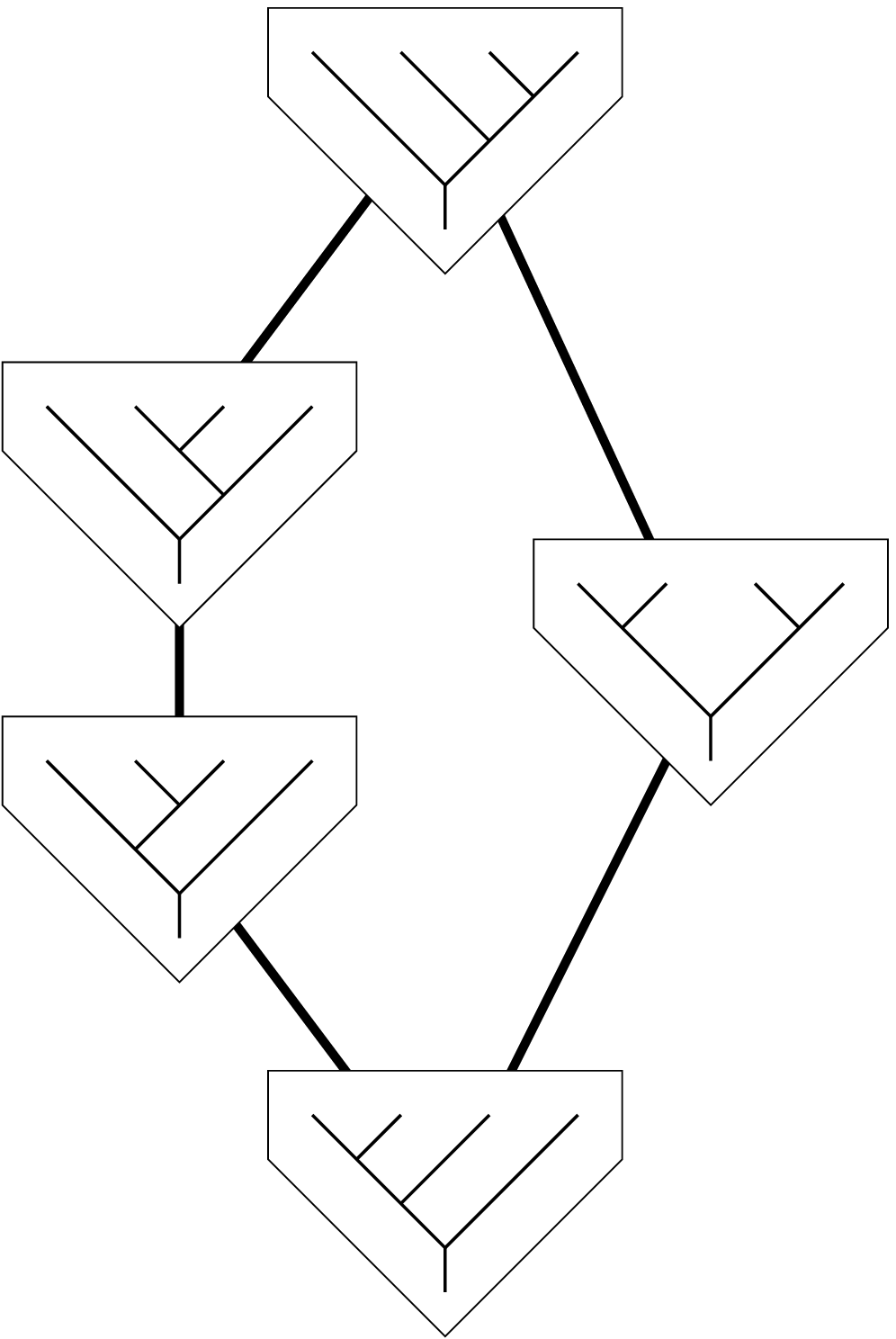,width=4cm} 
    \caption{The Tamari poset $\TT(2)$}
    \label{tam2}
  \end{center}
\end{figure}

This poset is called the Tamari poset of degree $n$, denoted by
$\TT(n)$. It is known to be a lattice. The lattice $\TT(2)$ is
depicted in Figure \ref{tam2}. 

The left-right symmetry of trees is an anti-automorphism of this
poset, sending the minimal element to the maximal element. 

The minimal element of $\TT(n)$ will be denoted by $\zer$ and the
maximal element by $\one$.

\begin{lemma}
  \label{slashlemma}
  For any $T^1,T^2$ in $\TT(n)$, the map $(s^1,s^2)\mapsto s^1
  \backslash s^2$ is a bijection from the product of the intervals
  $[T^1,\one] \times [T^2,\one]$ to the interval $[T^1 \backslash
  T^2,\one]$. 
\end{lemma}
\begin{proof}
  This is quite obvious from the definition of the partial order, as
  the covering relations preserve the fact that a tree can be
  written $s^1 \backslash s^2$. 
\end{proof}

\section{Dendriform algebras}

The notion of dendriform algebra was introduced by Loday, see
\cite{loday}. Let us recall the axioms. 

A dendriform algebra over some field $\kk$ is a vector space over
$\kk$ with two maps $\prec,\succ\, : \kk\otimes \kk \to \kk$
satisfying the following equations:
\begin{align}
  \label{axiom1}
  (x \prec  y) \prec z &= x \prec ( y \prec z)+x \prec  (y \succ z),\\
  \label{axiom2}
  x \succ ( y \prec z) &= (x \succ  y) \prec z,\\
  \label{axiom3}
   x \succ (y \succ z) &= (x \succ y) \succ z +(x \prec y) \succ z. 
\end{align}

These relations implies that the map $*$ defined by $x *y = x \prec
y+x \succ y$ is associative. 

There is a nice description of the free dendriform algebra on one
generator in terms of planar binary trees, see
\cite{loday,lodayronco1}. In particular, the underlying vector space
is $\kk \YY^+$. One can define the operations $\prec$ and $\succ$ on
$\kk \YY^+$. The product $*$ can be extended to $\kk \YY$ and has an
inductive definition as follows.

\begin{proposition}
  \label{easy}
  The tree $|$ is a unit for $*$. For all $T^1,T^2,T^3,T^4$ in $\YY$, one has
  \begin{equation}
    (T^1 _{Y} T^2)*(T^3 _{Y} T^4)
    =((T^1 _{Y} T^2)*T^3) _{Y} T^4+ T^1_{Y}(T^2*(T^3 _{Y} T^4)). 
  \end{equation}
\end{proposition}

There is also a simple expression for the product $*$ in $\kk \YY$
which uses the Tamari poset \cite[Eq. (2)]{lodayronco2}. 

\begin{proposition}
  Let $S$ and $T$ be in $\YY$. One has the following relation in $\kk \YY$:
  \begin{equation}
    S * T = \sum_{ S/T  \leq U \leq S\backslash T} U. 
  \end{equation}
\end{proposition}

We will need the following Lemma. 
\begin{lemma}
  \label{moblemma}
  For any $T^1,T^2$ in $\TT(n)$, the product of the intervals
  $[\zer,T^1]$ and $[\zer,T^2]$ is exactly the interval $[\zer,T^1
  \backslash T^2]$. 
\end{lemma}
For its proof, see for example \cite[Th. 29 \& 30]{hivnov}. 

\section{The Dendriform operad}

As a reference on operads and anticyclic operads, the reader may wish
to consult \cite{markl,book}. 

In this paper, we will only consider non-symmetric operads. A
non-symmetric operad $\PP$ in the category of vector spaces over $\kk$
is a collection of vector spaces $\PP(n)$ for $n\geq 1$, a collection
of maps $\circ_i : \PP(n)\otimes\PP(m) \to \PP(n+m-1)$ for $1\leq i
\leq n$ and a unit $1$, satisfying axioms modelled after the
composition of some multi-linear map at some place $i$ inside another
multi-linear map. The unit $1$ plays the rôle of the identity map in
the composition of multi-linear maps. 

A anticyclic non-symmetric operad $\PP$ is a non-symmetric operad
together with a linear map $\tau$ on each $\PP(n)$ such that
$\tau^{n+1}=\Id$ and the following relations hold for $a\in\PP(n)$ and
$b \in \PP(m)$:
\begin{align}
  \tau(1)&=-1,\\
  \tau(a \circ_n b)&=-\tau(b)\circ_1 \tau(a),\\
  \tau(a \circ_i b)&=\tau(a)\circ_{i+1} b \quad \text{ if }1\leq
  i<n. 
\end{align}

Let us now define the dendriform operad $\dend$. For all $n\geq 1$,
the space $\dend(n)$ is the vector space $\kk\YY(n)$ spanned by the
set of planar binary trees of degree $n$. The composition maps
$\circ_i$ can be described using shuffles of trees, see \cite[Prop
5.11]{loday}. The unit of the operad $\dend$ is the unique tree with
one vertex, denoted by $Y$. 

The operad $\dend$ is generated by two elements $\prec$ and $\succ$
with relations corresponding to Formulas
(\ref{axiom1},\ref{axiom2},\ref{axiom3}). These two elements should be
seen as the two elements of $\YY(2)$, namely $\prec$ is the tree
$\drt$ and $\succ$ is the tree $\gch$. 

Some of the combinatorial operations and products defined before can
be restated using the composition maps of the operad $\dend$. 

\begin{proposition}
  For all $T^1,T^2$ in $\YY^+$, one has the following relations:
  \begin{align}
    T^1 * T^2 &= ((\gch+\drt) \circ_2 T^2) \circ_1 T^1 ,\\
    T^1 \backslash T^2 &= T^1 \circ_{n_1} (\drt \circ_{2} T^2),\\ 
    T^1 / T^2 &= T^2 \circ_1 (\gch \circ_1 T^1). 
  \end{align}
  where $n_1$ is the degree of $T^1$. 
\end{proposition}

The following Theorem was proved in \cite[Thm. 4.1]{anticyclic} (in
some equivalent form). 

\begin{theorem}
  There exists a unique structure of anticyclic non-symmetric operad
  on $\dend$ such that 
  \begin{equation}
    \tau(\drt)=\gch\quad \text{and}\quad
    \tau(\gch)=-(\drt+\gch). 
  \end{equation}
\end{theorem}

The main aim of the present article is to gain some understanding of
the induced cyclic actions on $\dend(n)$. 

\section{Quivers}

\subsection{Quiver with relations from a poset}

Recall that a quiver $Q$ is a set of vertices $V$ and a set of arrows
$A$ with two maps from $A$ to $E$ giving the source and target of each
arrow. 

Then a module $M$ over $Q$ is a set of vector spaces $M_{v}$ for each
$v$ in $V$ and a set of maps $f_{v,w}$ from $M_{v}$ to $M_{w}$ for
each arrow in $A$ with source $v$ and target $w$. Modules over a
quiver $Q$ form an abelian category, denoted by $\mod(Q)$. 

One can restrict this category by imposing further conditions on the
composition of the maps $f_{v,w}$. For example, if $\Poset$ is a finite
poset, one can define a quiver $Q_\Poset$ with vertices the elements of $\Poset$
and arrows the covering relations of $\Poset$. That is to say, there is an
arrow from $v$ to $w$ in $Q_\Poset$ if and only if $v \leq w$ in $\Poset$ and
there is no element $u$ in $\Poset$ such that $v <u < w$. 

Then one can consider the category $\mod(\Poset)$ of modules over the quiver
$Q_\Poset$ such that for any pair $v \leq w$ in $\Poset$ and any two sequences
of arrows $v=u_0 \to u_1 \to u_2 \to \dots \to u_k=w$ , $v=u'_0 \to
u'_1 \to u'_2 \to \dots \to u'_\ell=w$ in $Q_\Poset$, one has the relation
\begin{equation}
  f_{u_0,u_1}f_{u_1,u_2} \dots f_{u_{k-1},u_k}
= f_{u'_0,u'_1}f_{u'_1,u'_2} \dots f_{u'_{\ell-1},u'_\ell},
\end{equation}
where composition of maps is denoted by concatenation. Then the
category $\mod(\Poset)$ is also an abelian category. As $\Poset$ is assumed
finite, this abelian category is known to have finite cohomological
dimension. 

\subsection{Derived category and Coxeter transformation}

Let $\Dmod(\Poset)$ be the bounded derived category of $\mod(\Poset)$. 

This derived category has a canonical self-equivalence which is called
the Auslander-Reiten translation, see \cite{happel, lenzing}. It is known that
this functor induces a map on the Grothendieck group $K_0$ of the
derived category. This map is called the Coxeter transformation. This
Grothendieck group has a natural basis indexed by the elements of $\Poset$,
corresponding to the images of simple modules of $\mod(\Poset)$ in the
derived category. 

We will denote by $\theta$ the Coxeter transformation in the
Grothendieck group of the derived category $\Dmod(\Poset)$. 

Let $L$ be the matrix defined by $L_{v,w}=1$ if and only if $v \leq w$
in $\Poset$. Then it is known that
\begin{proposition}
  \label{coxet}
  The matrix of the Coxeter transformation $\theta$ in the natural
  basis of $K_0$ is given by $-L ({L}^t)^{-1}$. 
\end{proposition}

Remark that $\theta$ is clearly an invertible map. 

From now on, this construction will be used for the Tamari posets
$\TT(n)$. In particular $\theta$ denotes the Coxeter transformation
for some Tamari poset $\TT(n)$, where $n$ should be clear from the
context. As the underlying set of $\TT(n)$ is $\YY(n)$, the action of
$\theta$ on $K_0(\Dmod(\TT(n)))$ can be
interpreted as an action on $\dend(n)$. 

\section{Periodicity Theorem}

\label{principale}

Here is the main result, relating the anticyclic structure of the
dendriform operad and the derived categories of modules on the Tamari
lattices. 

\begin{theorem}
  \label{main}
  On the vector space $\dend(n)$, one has the relation
  \begin{equation}
    \tau=(-1)^n \theta^2. 
  \end{equation}
\end{theorem}

The proof of this Theorem is done in the next section. Before this
proof, let us state a consequence. 

\begin{corollary}
  The Coxeter transformation $\theta$ in the Grothendieck group of the
  derived category $\Dmod(\TT(n))$ of modules on the Tamari lattice
  $\TT(n)$ satisfies $\theta^{2n+2}=\Id$. 
\end{corollary}

\begin{proof}
  As part of the anticyclic structure on $\dend$, it is known that
  $\tau^{n+1}=\Id$ on $\dend(n)$. 
\end{proof}

\subsection{Proof of the main theorem}

The strategy of proof is to find some inductive characterization of
the map $\tau$ and then to prove that the map $(-1)^n \theta^2$
satisfies the same induction. 

\begin{proposition}
  \label{char_tau}
  The collection of maps $\tau$ is uniquely defined by the following
  equations, for all $T, T^1$, $T^2$ in $\YY^+$. 
  \begin{align}
    \tau(Y)&=-Y,\\
    \tau(T^1 \backslash T^2)&=\tau(T^1) / \tau(T^2),\\
    \tau(T /Y )&=-Y * T. 
  \end{align}
\end{proposition}

\begin{proof}
  The fact that $\tau(Y)=-Y$ is by definition of an anticyclic operad. 

  Let us first prove that $\tau$ satisfies these equations, using the
  axioms of anticyclic operad and the known action of $\tau$ on $\gch$
  and $\drt$. One has
  \begin{equation*}
    \tau(T/Y)=\tau(\gch \circ_1 T)=\tau(\gch) \circ_2
    T=-(\gch+\drt)\circ_2 T=-Y * T. 
  \end{equation*}
  Let $n_1$ be the degree of $T^1$. One also has
  \begin{equation*}
    \tau(T^1\backslash T^2)=\tau(T^1 \circ_{n_1}(\drt \circ_2 T^2))
    =(\tau(T^2)\circ_1 \gch)\circ_1 \tau(T^1)=\tau(T^1)/\tau(T^2). 
  \end{equation*}

  The proof of uniqueness is an easy induction on degree. Any tree $T$
  in $\YY^+$ which is not $Y$ can either be written $T^1 \backslash
  T^2$ for some trees in $\YY^+$ of smaller degrees, or has the shape
  $T'/Y$ for some tree $T'$ of smaller degree. This allows to define
  $\tau$  by induction. 
\end{proof}

Let us now prove some properties of $\theta$ and deduce from them
properties of $\theta^2$. 

\begin{proposition}
  \label{char_theta}
  The collection of maps $\theta$ satisfy the following relations, for
  all $T^1$, $T^2$ in $\YY$. 
  \begin{align}
    \theta(|)&=-|,\\
    \theta(Y)&=-Y,\\
    \theta(T^1 \backslash T^2)&=-\theta(T^1) * \theta(T^2),\\
    \theta(T^1 * T^2)&=-\theta(T^1) / \theta(T^2),\\
    \theta^{-1}(T^1 / T^2)&=-\theta^{-1}(T^1) * \theta^{-1}(T^2),\\
    \theta^{-1}(T^1 * T^2)&=-\theta^{-1}(T^1) \backslash \theta^{-1}(T^2). 
  \end{align}
\end{proposition}

\begin{proof}
  It is clear that $\theta(|)=-|$ and $\theta(Y)=-Y$. The equations
  for $\theta^{-1}$ are obvious consequences of the equations for
  $\theta$. It is enough to prove one of the equations for $\theta$ as
  they are related by conjugation by the left-right symmetry of trees. 
  Let us prove the first one. By the definition of $-\theta$ from
  Proposition \ref{coxet}, it is the composite of the matrices $L$ and
  $({L}^t)^{-1}$. By Lemma \ref{slashlemma}, the action of ${L}^t$
  preserves the $\backslash$ product. Hence this is also true for its
  inverse. By Lemma \ref{moblemma}, the action of $L$ maps the
  $\backslash$ product to the $*$ product. Hence $-\theta$ maps the
  $\backslash$ product to the $*$ product. This proves the
  Proposition. 
\end{proof}

Remark that the conditions in Prop. \ref{char_theta} in fact uniquely
determine the collection of maps $\theta$. We will not need that fact. 

\begin{corollary}
  \label{coro1}
  For all $T^1,T^2$ in $\YY$ of degree $n_1,n_2$, one has the
  following relation
  \begin{equation}
    (-1)^n \theta^2(T^1 \backslash T^2)=(-1)^{n_1}
    \theta^2(T^1) / (-1)^{n_2}\theta^2(T^2),
  \end{equation}
  where $n=n_1+n_2$ is the degree of $T^1/T^2$. 
\end{corollary}

We need another property of $\theta$. 

\begin{proposition}
  For all $T$ in $\YY$ of degree $n$, one has
  \begin{equation}
    \theta(T/Y)=(-1)^n Y \backslash \theta^{-1} (T)
    \quad\text{and}\quad
    \theta^{-1}(Y\backslash T)=(-1)^n \theta(T)/Y. 
  \end{equation}
\end{proposition}

\begin{proof}
  The proof is by induction on the degree of $T$. It is enough to
  prove one of the equations as they are obviously equivalent. The
  Proposition is clearly true for small degrees. Assume that $T$ can
  be written $T^1_{Y} T^2$ with $T^1$ of degree $n_1$ and $T^2$ of
  degree $n_2$ in $\YY$ with $n_1+n_2+1=n$. Then one has $T*Y =
  T/Y+(T^1/Y)\backslash(T^2*Y)$. Hence one gets on the one hand
  \begin{equation}
    \theta(T/Y)=\theta(T*Y)-\theta((T^1/Y)\backslash(T^2*Y)). 
  \end{equation}
  Then using twice Proposition \ref{char_theta}, this becomes
  \begin{equation}
    \theta(T)/Y+\theta(T^1/Y)*\theta(T^2*Y). 
  \end{equation}
  Using again Proposition \ref{char_theta} and the fact that
  $T=(T^1/Y)\backslash T^2$, this is 
  \begin{equation}
        -(\theta(T^1/Y)*\theta(T^2))/Y+\theta(T^1/Y)*(\theta(T^2)/Y). 
  \end{equation}
  Then using the induction hypothesis on $T^1$, one gets
  \begin{equation}
    (-1)^{n_1+1}((Y\backslash\theta^{-1}(T^1))*\theta(T^2))/Y+(-1)^{n_1}(Y\backslash\theta^{-1}(T^1))*(\theta(T^2)/Y). 
  \end{equation}
  
  On the other hand, using the fact that $T=T^1/(Y\backslash T^2)$ and
  Proposition \ref{char_theta}, one has
  \begin{equation}
    (-1)^n Y \backslash \theta^{-1} (T)
    =(-1)^{n_1+n_2} Y \backslash
    (\theta^{-1}(T^1)*\theta^{-1}(Y\backslash T^2)). 
  \end{equation}
  Using the induction hypothesis for $T^2$, one gets
  \begin{equation}
    =(-1)^{n_1} Y \backslash (\theta^{-1}(T^1)*(\theta(T^2)/Y)). 
  \end{equation}
  Then using Proposition \ref{easy} for $a=\theta^{-1}(T^1)$ and
  $b=\theta(T^2)$, the induction step is done. 
\end{proof}

\begin{corollary}
  \label{coro2}
  For all $T$ in $\YY$ of degree $n$, one has
  \begin{equation}
    (-1)^{n+1} \theta^2 (T / Y)=-Y * T. 
  \end{equation}
\end{corollary}

From Corollaries \ref{coro1} and \ref{coro2} and by Proposition
\ref{char_tau}, one gets a proof of Theorem \ref{main}. 

\bibliographystyle{plain}
\bibliography{tamari}

\end{document}